\documentclass[twoside,12pt]{article}
\usepackage{color}

\usepackage{graphics}
\usepackage{graphicx}

\usepackage{amssymb,amsfonts}
\usepackage{amsmath}
\usepackage{amsthm}

\def\binomh#1#2{ \scalebox{.3}[1.2]{\textbf{)}}{\genfrac{}{}{0pt}{}{#1}{#2}}\scalebox{.3}[1.2]{\textbf{(}} }

\def\binomhh#1#2{ \scalebox{.4}[1.7]{\textbf{)}}{\genfrac{}{}{0pt}{}{#1}{#2}}\scalebox{.4}[1.7]{\textbf{(}} }

\def\ve{\varepsilon}

\newtheorem{theorem}{Theorem}
\newtheorem{lemma}{Lemma}[section]
\newtheorem{cor}{Corollary}
\newtheorem{rem}{Remark}

\title{\bf Alternating sums in hyperbolic Pascal triangles 
}

\author{L\'aszl\'o N\'emeth\footnote{University of West Hungary,  Institute of Mathematics, Hungary. \textit{nemeth.laszlo@nyme.hu}}, 
 L\'aszl\'o Szalay\footnote{Department of Mathematics and Informatics, J. Selye University, Hradna ul. 21., 94501 Komarno, Slovakia} \footnote{University of West Hungary,  Institute of Mathematics, Hungary.  \textit{szalay.laszlo@nyme.hu}}}
\date{}

\begin{document}

\maketitle

\begin{abstract}
A new generalization of Pascal's triangle, the so-called hyperbolic Pascal triangles were introduced in \cite{BNSz}. The mathematical background goes back to the regular mosaics in the hyperbolic plane. The alternating sum of elements in the rows was given in the special case $\{4,5\}$ of the hyperbolic Pascal triangles. In this article, we determine the alternating sum generally in the hyperbolic Pascal triangle corresponding to $\{4,q\}$ with $q\ge5$. \\[1mm]
\end{abstract}

\section{Introduction}\label{sec:introduction}
 
In the hyperbolic plane there are infinite types of regular mosaics (see, for example \cite{C}), they are denoted by Schl\"afli's symbol $\{p,q\}$, where $(p-2)(q-2)>4$. Each regular mosaic induces a so called hyperbolic Pascal triangle (see \cite{BNSz}), following and generalizing the connection between the classical Pascal's triangle and the Euclidean regular square mosaic $\{4,4\}$.
For more details see \cite{BNSz}, but here we also collect some necessary information. 

\begin{figure}[h!]
 \centering
 \includegraphics[width=0.99\linewidth]{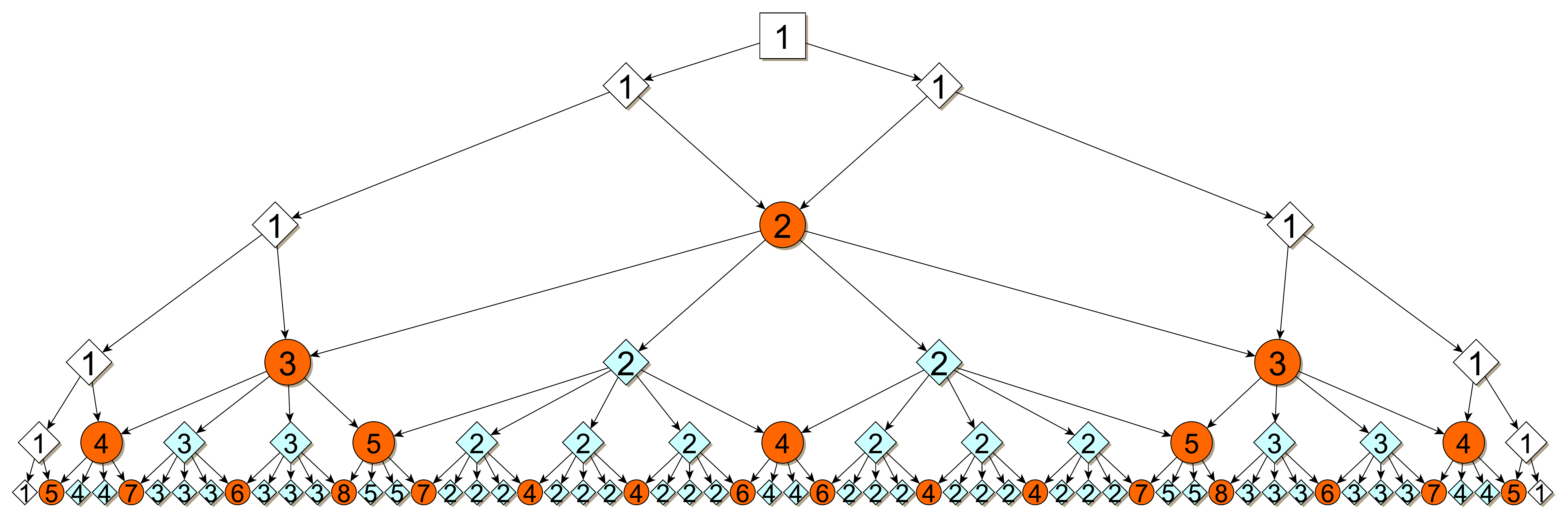}
 \caption{Hyperbolic Pascal triangle linked to $\{4,6\}$ up to row 5}
 \label{fig:Pascal_46_layer5}
\end{figure}

There are several approaches to generalize the Pascal's arithmetic triangle (see, for instance \cite{BSz}).
The hyperbolic Pascal triangle based on the mosaic $\{p,q\}$ can be figured as a digraph, where the vertices and the edges are the vertices and the edges of a well defined part of the lattice $\{p,q\}$, respectively, further the vertices possesses a value giving the number of different shortest paths from the base vertex. Figure~\ref{fig:Pascal_46_layer5} illustrates the hyperbolic Pascal triangle when $\{p,q\}=\{4,6\}$. 
Generally, for $\{4,q\}$ the base vertex has two edges, the leftmost and the rightmost vertices have three, the others have $q$ edges. The square shaped cells surrounded by appropriate edges are corresponding to the regular squares in the mosaic.
Apart from the winger elements, certain vertices (called ``Type A'' for convenience) have two ascendants and $q-2$ descendants, the others (``Type B'') have one ascendant and $q-1$ descendants. In the figures we denote the vertices type $A$ by red circles and the vertices type $B$ by cyan diamonds, further the wingers by white diamonds. The vertices which are $n$-edge-long far from the base vertex are in row $n$. 

The general method of drawing is the following. Going along the vertices of the $j^{th}$ row, according to type of the elements (winger, $A$, $B$), we draw appropriate number of edges downwards (2, $q-2$, $q-1$, respectively). Neighbour edges of two neighbour vertices of the $j^{th}$ row meet in the $(j+1)^{th}$ row, constructing a vertex with type $A$. The other descendants of row $j$ in row $j+1$ have type $B$.
In the sequel, $\binomh{n}{k}$ denotes the $k^\text{th}$ element in row $n$, which is either the sum of the values of its two ascendants or the value of its unique ascendant. We note, that the hyperbolic Pascal triangle has the property of vertical symmetry. 

It is well-known that the alternating sum $\sum_{i=0}^{n}(-1)^i\,{n\choose i}$ of row $n$ in the classical Pascal's triangle is zero ($n\geq1$). In \cite{BNSz} we showed that the alternating sum $\sum_{i}(-1)^i\,\binomh{n}{i}$ is either 0 (if $n\equiv1\;(\bmod\ 3)$) or 2 (otherwise, with $n\ge5$) for case $\{4,5\}$. In this paper, we determine an explicit form for the alternating sums generally for $\{4,q\}$ $(q\geq5)$. If one considers the  result with $q=4$, it returns 0 according to the classical Pascal's triangle. The definitions, the signs, the figures, and the method strictly follow the article \cite{BNSz}. 

\section{Main Theorems}
From this point we consider the hyperbolic Pascal triangle based on the mosaic $\{4,q\}$ with $q\geq5$. 
Denote by $s_n$ and $\hat{s}_n$ the number of the vertices, and the sum of the elements in row $n$, respectively. 

The sequences $\{s_n\}$ and $\{\hat{s}_n\}$ can be given (see again \cite{BNSz}) by the  ternary homogenous recurrence relations 
\begin{equation}\label{eq:recsn}
s_n=(q-1)s_{n-1}-(q-1)s_{n-2}+s_{n-3}\qquad (n\ge4),
\end{equation}
(the initial values are $ s_1=2,\;s_2=3,\;s_3=q$) and 
\begin{equation*}
\hat{s}_n=q\hat{s}_{n-1}-(q+1)\hat{s}_{n-2}+2\hat{s}_{n-3}\qquad (n\ge4),
\end{equation*}
(the initial values are $\hat{s}_1=2,\;\hat{s}_2=4,\;\hat{s}_3=2q$), respectively. 
  
Let $\widetilde{s}_{n}$  be the alternating sum of elements of the hyperbolic Pascal triangle (starting with positive coefficient) in row $n$, and we distinguish the even and odd cases. 

\begin{theorem} \label{th:even}
Let $q$ be even. Then
$$
\widetilde{s}_n=\sum_{i=0}^{s_n-1}(-1)^i\,\binomhh{n}{i}=
\left\{
\begin{array}{clll}
0, & {\rm if} & n=2t+1, & n\ge1,\\
-2(5-q)^{t-1}+2, & {\rm if} & n=2t, & n\ge2,
\end{array}
\right.
$$
hold, further  $\widetilde{s}_0=1$. 
\end{theorem}


\begin{cor}
In case of $q=6$ we deduce 
$$
\widetilde{s}_n=
\left\{
\begin{array}{clll}
0, & {\rm if} & n\ne 4t, & n\ge1,\\
4, & {\rm if} & n=4t, & n\ge4.
\end{array}
\right.
$$ 
\end{cor}

\begin{rem}
For $q=4$ Theorem~\ref{th:even} would return with 0, providing the known result $ \widetilde{s}_n=0$ for the original Pascal's triangle.  
  
\end{rem}

\begin{theorem} \label{th:odd}
Let $q\ge5$ be odd. Then $\widetilde{s}_0=1$, further
$$
\widetilde{s}_n=\sum_{i=0}^{s_n-1}(-1)^i\,\binomhh{n}{i}=
\left\{
\begin{array}{clll}
0, & {\rm if} & n=3t+1,& n\ge1 ,\\
\phantom{2}(-2)^t(q-5)^{t-1}+2, & {\rm if} & n=3t-1,& n\ge n_1,\\
2(-2)^t(q-5)^{t-1}+2, & {\rm if} & n=3t,& n\ge n_2,
\end{array}
\right.
$$
where $(n_1,n_2)=(2,3)$ and $(5,6)$ if $n>5$ and $n=5$, respectively. In the latter case $\widetilde{s}_2=0$, $\widetilde{s}_3=-2$. 
$$
$$
\end{theorem}

We note, that by the help of $\hat{s}_n$ and $\widetilde{s}_n$ we can easily determine the alternate sum with the arbitrary weights $v$ and $w$. 
\begin{cor}
\begin{eqnarray*}
\widetilde{s}_{(v,w),n}&=&\sum_{i=0}^{s_n-1}\left(v{\delta_{0, \,i\bmod2}}+w{\delta_{1,\, i\bmod2}}\right)\binomhh{n}{i}=\frac{\hat{s}_n+\widetilde{s}_n}{2}v+\frac{\hat{s}_n-\widetilde{s}_n}{2}w\\
&=&
\frac{v+w}{2}\hat{s}_n+\frac{v-w}{2}\widetilde{s}_n,
\end{eqnarray*}
where $\delta_{j,i}$ is the Kronecker delta.
\end{cor}

\section{Proofs of Theorems}

Since the hyperbolic Pascal triangle has a symmetry axis, if $s_n$ is even then the alternating sum is zero. In the case when $s_n$ is odd (in the sequel, we assume it), the base of both proofs is to consider the vertices type $A$ and $B$ of row $n$ and to observe their influence on $\widetilde{s}_{n+2}$ or $\widetilde{s}_{n+3}$. We separate the contribution of each $\binomh{n}{i}$ individually, and then take their superposition. 

Let  $\widetilde{s}_n^{(A)}$ and $\widetilde{s}_n^{(B)}$ be the subsum of $\widetilde{s}_n$ restricted only to the elements of type $A$ and $B$, respectively. 
As $\binomh{n}{0}=\binomh{n}{s_n-1}=1$, then 
\begin{equation}\label{eq:snAB}
\widetilde{s}_{n}=\widetilde{s}_n^{(A)}+\widetilde{s}_n^{(B)}+2.
\end{equation}  

Using the notations of \cite{BNSz}, $x_A$ and $x_B$ denote the value of an element of type $A$ and $B$, respectively. (In the figures, we indicate them shortly by $x$.)
Their contributions to $\widetilde{s}_{n+k}$ $(k\geq1)$ are denote by ${\cal H}_k(x_A)$ and ${\cal H}_k(x_B)$, respectively, and for example, 
${\cal H}_k^{(A)}(x_A)$ and ${\cal H}_k^{(B)}(x_A)$ for the contribution of  $x_A$ from the row $n$ to the alternate sum of the row $n+k$ restricted to the elements of type $A$ and $B$, respectively. Similarly, ${\cal H}_k(1)$ shows the contribution of a winger element of row $n$ (having value $1$).
According to \cite{BNSz} we have 
\begin{eqnarray*}\label{eq:inf}
{\cal H}_k(x_A)&=&{\cal H}_k^{(A)}(x_A)+{\cal H}_k^{(B)}(x_A),\\
{\cal H}_k(x_B)&=&{\cal H}_k^{(A)}(x_B)+{\cal H}_k^{(B)}(x_B),\\
{\cal H}_k(1)&=&{\cal H}_k^{(A)}(1)+{\cal H}_k^{(B)}(1)+1.
\end{eqnarray*}

Clearly,  it is obvious that $\widetilde{s}_{0}=1$, $\widetilde{s}_{2}=0$ hold for all $q\geq5$.

\subsection{Proof for even $q$}

If $n=2t+1$ $(n\ge1)$, then in accordance with relation \eqref{eq:recsn} $s_n$ is even, otherwise odd. 
So, we may suppose that $n$ is also even.

\begin{figure}[h!]
 \centering
 \includegraphics[width=0.99\linewidth]{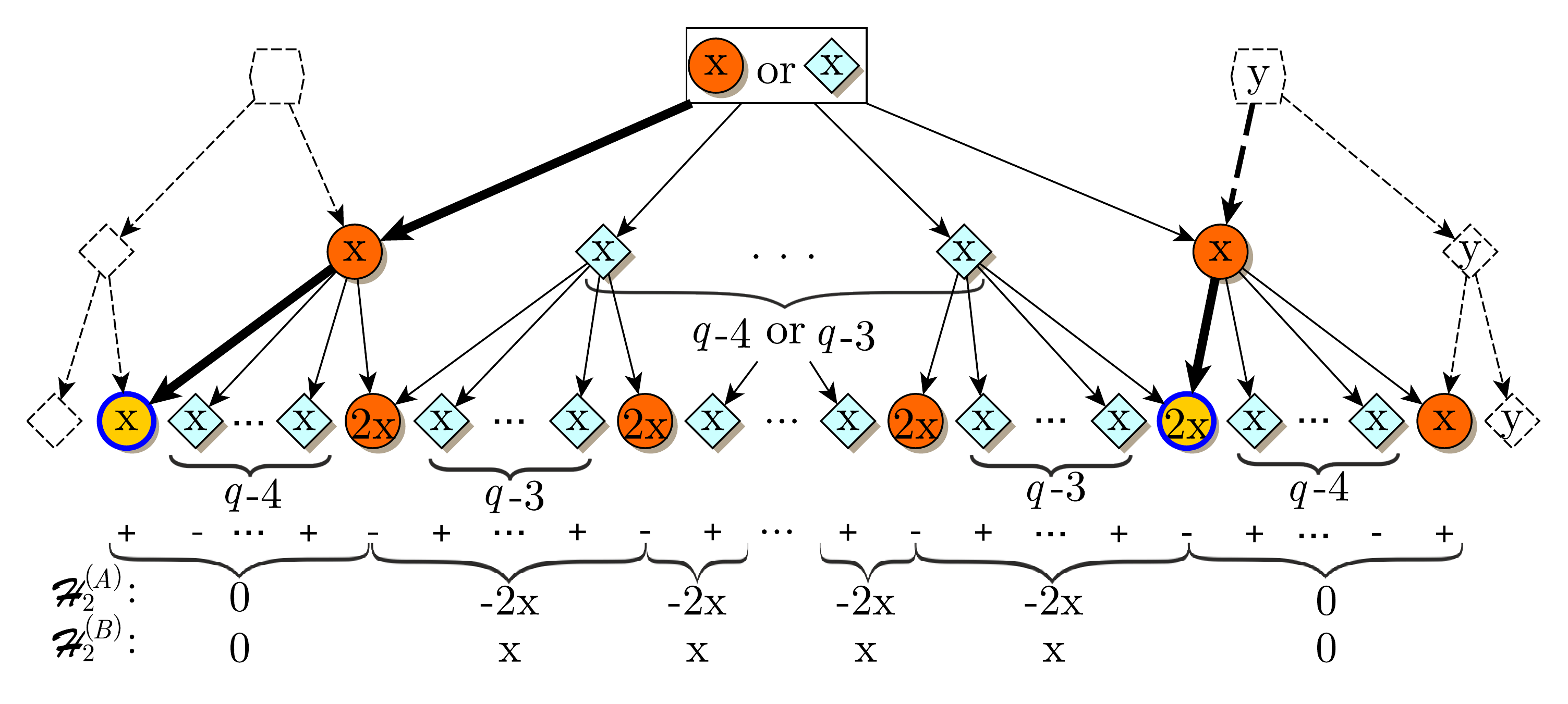}
\caption{The influence ${\cal H}_2(x_A)$ and ${\cal H}_2(x_B)$}
 \label{fig:4q_even_HAB2}
\end{figure}

Figure~\ref{fig:4q_even_HAB2} shows the contributions of  
$x_A$ and $x_B$ from the row $n$ to the alternating sum of the row $n+2$. From the growing method of the hyperbolic Pascal triangle, a vertex type $A$ in row $i$ generates  $q-4$ vertices type $B$ in row $i+1$, if the vertex is type $B$ then it has $q-3$ generated vertices type $B$ in row $i+1$.  The value of a vertex is either the value of its ascendant or the sum of values of its two ascendants. Then we consider the value $2x_A$ (and $2x_B$) of vertices type $A$ in rows $n+2$  as $x_A+x_A$ (and $x_B+x_B$). In Figure~\ref{fig:4q_even_HAB2} we drew the values of the alternating sums belonging to $x$ (in row $n+2$) in blocks. The last but one row shows the values ${\cal H}_2^{(A)}(x_A)$ (or ${\cal H}_2^{(A)}(x_B)$), the last row shows  ${\cal H}_2^{(B)}(x_A)$ (or ${\cal H}_2^{(B)}(x_B)$) in blocks.

Put $\ve=\pm1$ and  $\delta=\pm1$ in case of the vertex type $A$ and $B$ is being considered, respectively. (In Figure~\ref{fig:4q_even_HAB2}, they are $+1$, it is the first sign in row $n+2$.)  Now we obtain, by observing  Figure~\ref{fig:4q_even_HAB2}, that  
\begin{eqnarray*}
{\cal H}_2^{(A)}(x_A) &=& \ve\cdot(-2x(q-4))= -2(q-4)\ve  x_A,\\
{\cal H}_2^{(B)}(x_A) &=& \ve\cdot x(q-4)= (q-4)\ve  x_A,\end{eqnarray*}
and
\begin{eqnarray*}
{\cal H}_2^{(A)}(x_B) &=& \delta\cdot (-2x(q-3))= -2(q-3)\delta x_B,\\
{\cal H}_2^{(B)}(x_B) &=& \delta\cdot x(q-3)  =(q-3)\delta x_B.
\end{eqnarray*}

\begin{figure}[h!]
 \centering
 \includegraphics[width=5cm]{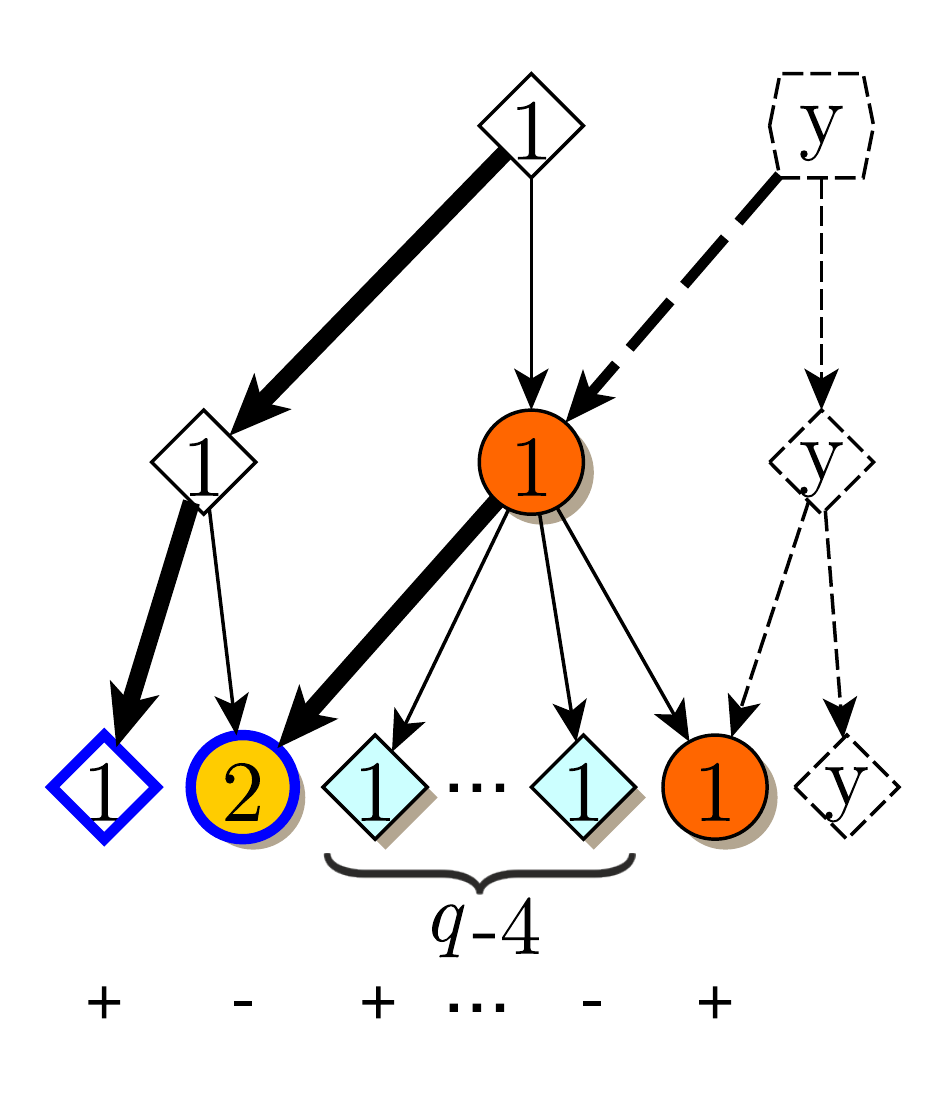}
\caption{The influence ${\cal H}_2(1)$}
 \label{fig:4q_even_H12}
\end{figure}

The Figure~\ref{fig:4q_even_H12} shows the contributions of the influence of left winger element of row $n$ to row $n+2$.
For the right winger elements the situation is the same, thanks to the vertical symmetry. Thus  
\begin{eqnarray*}
{\cal H}_2^{(A)}(1) &=& -1,\\
{\cal H}_2^{(B)}(1) &=& 0.
\end{eqnarray*}

We have given the influence of an element located in row $n$ on row $n+2$. Let suppose that $y$ is the value of the neighbour element of $x$ in row $n$. Clearly, $y$ has also influence on row $n+1$. The signs of $x$ and $y$ in the alternating sum in row $n$ are different and the signs of left hand side of their influence structures are also different in row $n+2$. In the figures, the leftmost element of influence structures are highlighted. (The  signs of rightmost values of influence structures are also different.) The situation is the same in case of the winger elements. Thus, according to \cite{BNSz} we can give the changing of the alternating sums from row $n$ to row $n+2$.        

Summarising the results,  we obtain the system of recurrence equations   
\begin{eqnarray}
\widetilde{s}_{n+2}^{(A)}&=&-2(q-4)\widetilde{s}_n^{(A)}-2(q-3)\widetilde{s}_n^{(B)}-2,\qquad (n\ge0),\label{eq:even_sA}\\
\widetilde{s}_{n+2}^{(B)}&=&\phantom{-2}(q-4)\widetilde{s}_n^{(A)}+\phantom{2}(q-3)\widetilde{s}_n^{(B)},\phantom{+1}\,\, \qquad (n\ge0).\label{eq:even_sB}
\end{eqnarray}

Now we apply the following lemma (see \cite{BNSz}).

\begin{lemma}\label{lemma:2seq}
Let $x_0$, $y_0$, further $u_i$, $v_i$ and $w_i$ ($i=1,2$) be complex numbers such that $a_2b_1\ne0$. Assume that
the for $n\ge n_0$ terms of the sequences $\{x_n\}$ and $\{y_n\}$ satisfy
\begin{eqnarray*}
x_{n+1}&=&u_1x_n+v_1y_n+w_1,\\
y_{n+1}&=&u_2x_n+v_2y_n+w_2.
\end{eqnarray*}
Then for both sequences
\begin{equation*}\label{nomix}
z_{n+3}=(u_1+v_2+1)z_{n+2}+(-u_1v_2+u_2v_1-u_1-v_2)z_{n+1}+(u_1v_2-u_2v_1)z_n
\end{equation*}
holds ($n\ge n_0$).
\end{lemma}

Thus we obtain
$\widetilde{s}_{n+6}^{(A)}=(6-q)\widetilde{s}_{n+4}^{(A)}+(q-5)\widetilde{s}_{n+2}^{(A)}$ and  $\widetilde{s}_{n+6}^{(B)}=(6-q)\widetilde{s}_{n+4}^{(B)}+(q-5)\widetilde{s}_{n+2}^{(B)}$. 
Obviously, $\widetilde{s}_{2}^{(A)}=-2$ and $\widetilde{s}_{2}^{(B)}=0$ fulfil. Using \eqref{eq:even_sA} and \textbf{\eqref{eq:even_sB}} we gain $\widetilde{s}_{4}^{(A)}=4q-18$ and  $\widetilde{s}_{4}^{(B)}=-2(q-4)$.

From \eqref{eq:snAB} we conclude
\begin{equation}\label{eq:even_recu}
\widetilde{s}_{n+6}=(6-q)\widetilde{s}_{n+4}+(q-5)\widetilde{s}_{n+2},\qquad (n\ge0),
\end{equation}
where $\widetilde{s}_{0}=1$, $\widetilde{s}_{2}=0$ and $\widetilde{s}_{4}=2(q-4)$.

The characteristic equation of \eqref{eq:even_recu} is
\begin{equation*}
\widetilde{p}(x)=x^4+(q-6)x^2-(q-5)=\left(x^2-(5-q)\right)(x^2-1).
\end{equation*}
Further, we have
$\widetilde{s}_{2t}=\alpha(5-q)^t+\beta$, and from $\widetilde{s}_2$ and $\widetilde{s}_4$ we realize $\alpha=-2/(5-q)$, $\beta=2$, and then $\widetilde{s}_{2t}=(-2/(5-q))(5-q)^t+2=-2(5-q)^{t-1}+2$.

\begin{rem}
For all $n\geq1$, $\widetilde{s}_n = -\widetilde{s}_{n}^{(B)}$\!, because  $\widetilde{s}_i = -\widetilde{s}_{i}^{(B)}=0$ for $i=1,2,3$ and $\widetilde{s}_4 = -\widetilde{s}_{4}^{(B)}$\!.
\end{rem}

\subsection{Proof for odd $q$}

We examine rows  $n\ne3t+1$ $(n\geq2)$, because the number of element in row $n=3t+1$ is even (see relation \eqref{eq:recsn}).

\begin{figure}[h!]
 \centering
 \includegraphics[width=0.80\linewidth]{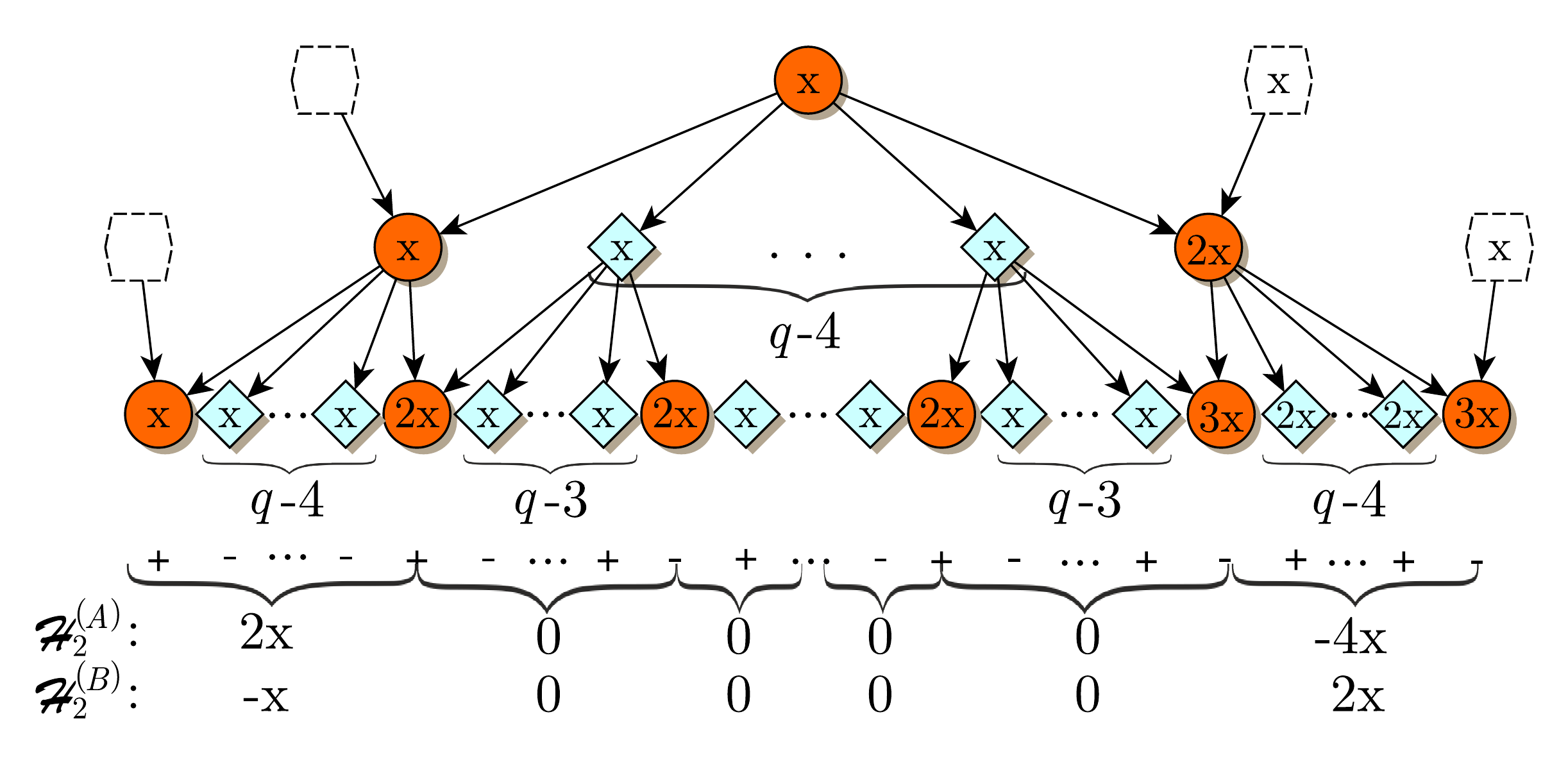}
\caption{The influence ${\cal H}_2(x_A)$}
 \label{fig:4q_odd_HA2}
\end{figure}

\begin{figure}[h!]
 \centering
 \includegraphics[width=0.80\linewidth]{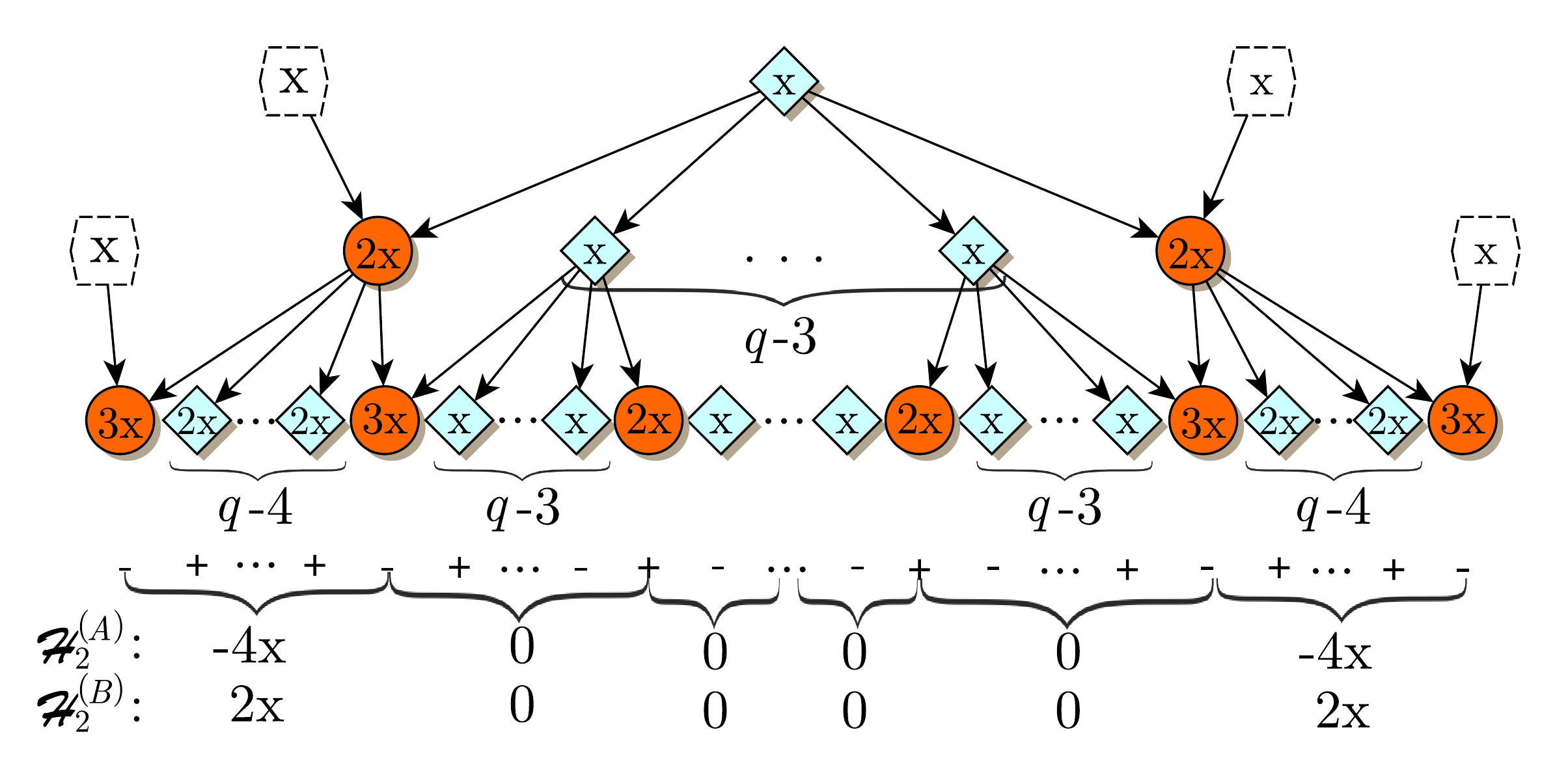}
\caption{The influence ${\cal H}_2(x_B)$}
 \label{fig:4q_odd_HB2}
\end{figure}

\begin{figure}[h!]
 \centering
 \includegraphics[width=0.97\linewidth]{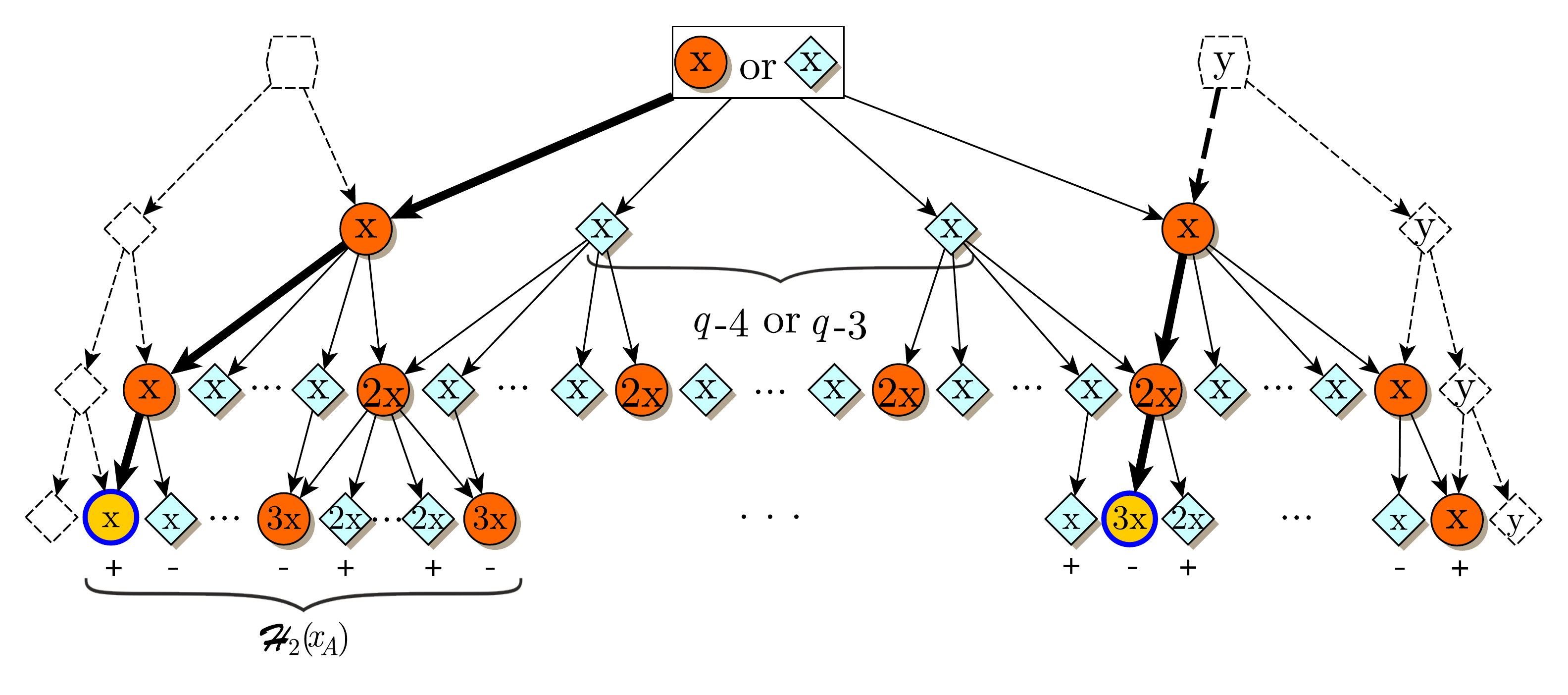}
\caption{The influence ${\cal H}_3(x_A)$ and ${\cal H}_3(x_B)$}
 \label{fig:4q_odd_HAB3}
\end{figure}

Here, apart from some details, we copy the treatment of the previous case. The first difference is that now we have to examine the influence of the elements from row $n$ on row $n+3$ because the nice property about the signs first appears three rows later. The influence structures ${\cal H}_3(x_A)$ and ${\cal H}_3(x_B)$ are rather complicated, so we split them into smaller parts. First we draw the structures ${\cal H}_2(x_A)$ and  ${\cal H}_2(x_B)$ when $x_A$ and $x_B$ are in row $n+1$, and we describe the influence of them on row $n+3$. Then we combine the result with the branches of ${\cal H}_3(x_A)$ and ${\cal H}_3(x_B)$. 
In Figures~\ref{fig:4q_odd_HA2} and \ref{fig:4q_odd_HB2} only ${\cal H}_2(x_A)$ and  ${\cal H}_2(x_B)$ can be seen, later in Figure~\ref{fig:4q_odd_HAB3} we consider the "skeleton" of ${\cal H}_3(x_A)$ and ${\cal H}_3(x_B)$.

\begin{figure}[h!]
 \centering
 \includegraphics[width=0.78\linewidth]{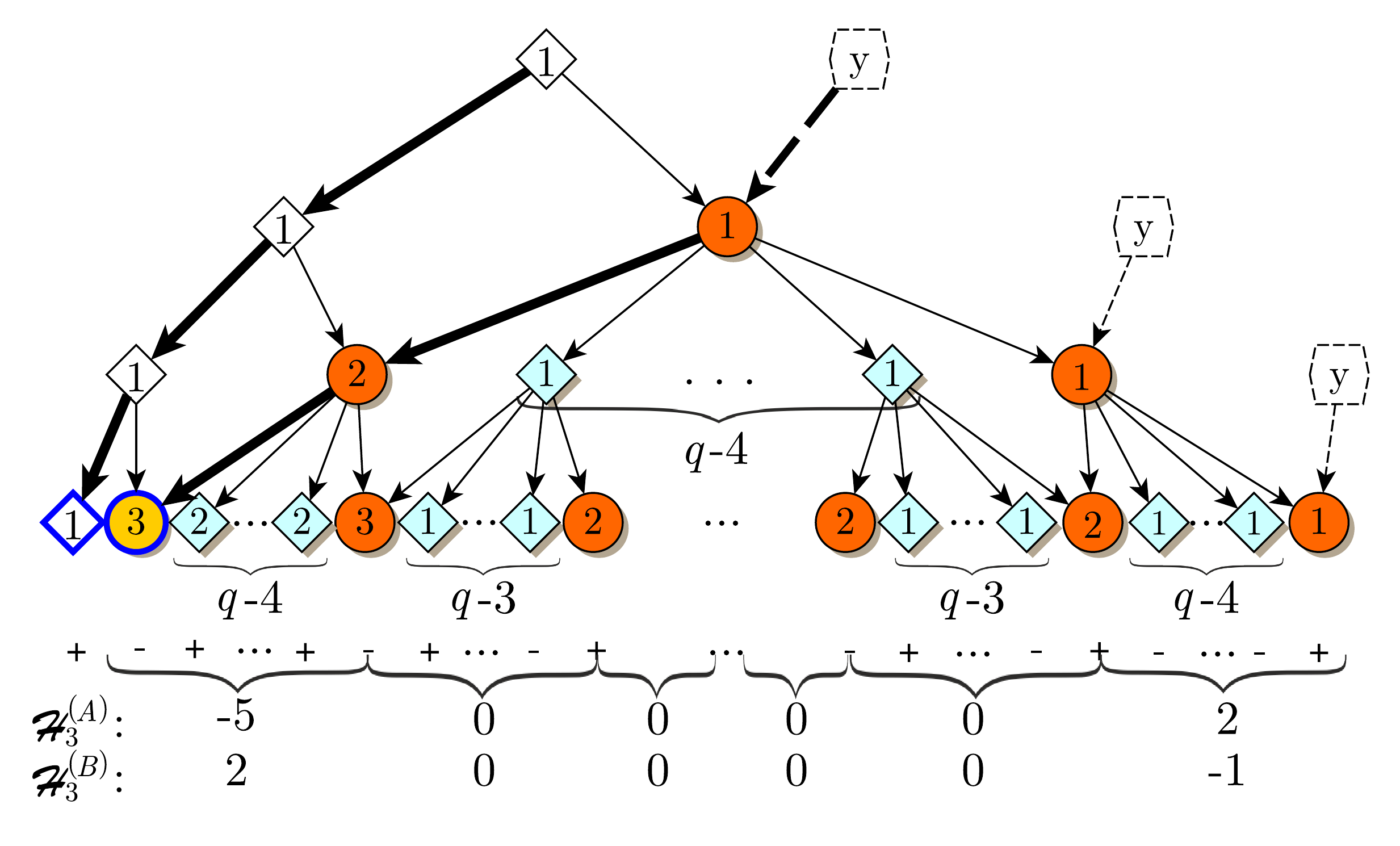}
\caption{The influence ${\cal H}_3(1)$}
 \label{fig:4q_odd_H13}
\end{figure}

From the figures one can derive the observations   
\begin{eqnarray*}
{\cal H}_3^{(A)}(x_A) &=& \ve\cdot(2x-4x-4x(q-4)+2x)=-4(q-4)\ve  x_A,\\
{\cal H}_3^{(B)}(x_A) &=& \ve\cdot(-x+2x+2x(q-4)-x)=2(q-4)\ve  x_A,
\end{eqnarray*}
\noindent and
\begin{eqnarray*}
{\cal H}_3^{(A)}(x_B) &=& \delta\cdot(2x-4x-4x(q-3)+2x)=-4(q-3)\delta x_B,\\
{\cal H}_3^{(B)}(x_B) &=& \delta\cdot(-x+2x+2x(q-3)-x)=2(q-3)\delta x_B.
\end{eqnarray*}
We also have
\begin{eqnarray*}
{\cal H}_3^{(A)}(1) &=& -5+2=-3,\\
{\cal H}_3^{(B)}(1) &=& 2-1=1.
\end{eqnarray*}

Combining the informations, it results the system of recurrence relations
\begin{eqnarray}\label{eq:_odd_recu01}
\widetilde{s}_{n+3}^{(A)}&=&-4(q-4)\widetilde{s}_n^{(A)}-4(q-3)\widetilde{s}_n^{(B)}+2\cdot(-3),\qquad (n\ge0),\\
\widetilde{s}_{n+3}^{(B)}&=&\;\;\;2(q-4)\widetilde{s}_n^{(A)}+2(q-3)\widetilde{s}_n^{(B)}+2\cdot1,\phantom{(-)}\qquad (n\ge0).\label{eq:_odd_recu02}
\end{eqnarray}

Lemma \ref{lemma:2seq} yields
$\widetilde{s}_{n+9}^{(A)}=(11-2q)\widetilde{s}_{n+6}^{(A)}+(2q-10)\widetilde{s}_{n+3}^{(A)}$ and  $\widetilde{s}_{n+9}^{(B)}=(11-2q)\widetilde{s}_{n+6}^{(B)}+(2q-10)\widetilde{s}_{n+3}^{(B)}$. 
Apparently, $\widetilde{s}_{2}^{(A)}=-2$, $\widetilde{s}_{2}^{(B)}=0$, and $\widetilde{s}_{3}^{(A)}=-6$, $\widetilde{s}_{3}^{(B)}=2$.
From the system \eqref{eq:_odd_recu01} and \eqref{eq:_odd_recu02} we gain
$\widetilde{s}_{5}^{(A)}=8q-38$,  $\widetilde{s}_{5}^{(B)}=-4q+18$ and
$\widetilde{s}_{6}^{(A)}=16q-78$,  $\widetilde{s}_{6}^{(B)}=-8q+38$.

By \eqref{eq:snAB} we realize
\begin{equation}\label{eq:odd_recu}
\widetilde{s}_{n+9}=(11-2q)\widetilde{s}_{n+6}+(2q-10)\widetilde{s}_{n+3},\qquad (n\ge0),
\end{equation}
where  $\widetilde{s}_{2}=0$, $\widetilde{s}_{3}=-2$ and $\widetilde{s}_{5}=4q-18$, $\widetilde{s}_{6}=8q-38$.

In case of $q=5$, the equation \eqref{eq:odd_recu} gives $\widetilde{s}_{n+9}=\widetilde{s}_{n+6}$ and $\widetilde{s}_{5}=2$, $\widetilde{s}_{6}=2$. So $\widetilde{s}_{n}=2$ holds if $n\geq 5$.
In case $q>5$, the characteristic equation of \eqref{eq:odd_recu} is
\begin{equation*}
\widetilde{p}(x)=x^6+(2q-11)x^3-(2q-10)=\left(x^3-2(5-q)\right)(x^3-1).
\end{equation*}

Thus 
$\widetilde{s}_{3t}$ or  $\widetilde{s}_{3t-1}=\alpha (2(5-q))^t+\beta$ fulfils for some $\alpha$ and $\beta$.
Finally, if $n=3t-1$, then from $\widetilde{s}_2$ and $\widetilde{s}_5$ we obtain $\alpha=1/(q-5)$ and $\beta=2$. Hence $\widetilde{s}_{3t-1}=(10-2q)^t/(q-5)+2$.
Otherwise, if $n=3t$, then from $\widetilde{s}_{3}$ and $\widetilde{s}_{6}$ we deduce $\alpha=2/(q-5)$ and $\beta=2$. Thus $\widetilde{s}_{3t}=2(10-2q)^t/(q-5)+2$.



\begin{thebibliography}{2015}

\bibitem{BNSz} 
H. Belbachir, L. N\'emeth and  L. Szalay, {Hyperbolic Pascal triangles}, \textit{Appl. Math. Comp.} {\bf 273} (2016), 453-464.

\bibitem{BSz} Belbachir, H.~and Szalay, L., On the arithmetic triangles, \u{S}iauliai Math.~Sem., {\bf 9} (17) (2014), 15-26. 

\bibitem{C} Coxeter, H.~S.~M., Regular honeycombs in hyperbolic space, Proc.~Int.~Congress of Math., Amsterdam, Vol. III. (1954), 155-169.

\end{thebibliography}
\end{document}